\documentclass[18pt]{article}
\usepackage{graphics}
\usepackage{amsmath}
\usepackage{amsfonts}
\usepackage{dsfont}
\usepackage{algorithm}
\input{pdfcolor}
\usepackage[noend]{algpseudocode}
\usepackage{adjustbox}
\usepackage{hyperref}

\usepackage{moresize}

\begin{document}

\title{Negative Latin Square Type Partial Difference Sets in Nonabelian Groups of Order 64}

\author{Andrew C. Brady\\Department of Mathematics and Computer Science\\University of Richmond\\Richmond, VA 23173\\email: {\tt andrew.brady1@richmond.edu}}

\date{}

\maketitle
\begin{abstract}

\noindent There exist few examples of negative Latin square type partial difference sets (NLST PDSs) in nonabelian groups. We present a list of 176 inequivalent NLST PDSs in 48 nonisomorphic, nonabelian groups of order 64. These NLST PDSs form 8 nonisomorphic strongly regular graphs. These PDSs were constructed using a combination of theoretical techniques and computer search, both of which are described. The search was run exhaustively on 212/267 nonisomorphic groups of order 64.

\end{abstract} 


\begin{section}{Introduction}

Let $G$ be a group, then we say that $D \subset G$ is a $(v,k,\lambda,\mu)$ partial difference set (PDS) if $v=|G|$, $k=|D|$, every element of $D$ appears in the multiset $\{d_1 d_2^{-1} | d_1, d_2 \in D\}$ $\lambda$ times, every nonidentity element of $G-D$ appears $\mu$ times, and the identity of $G$ appears $k$ times. Partial difference sets are important to understand better because they are equivalent to strongly regular graphs (SRGs) with regular automorphism groups, can be used to create rank 3 association schemes, and can be used to create projective 2-weight error correcting codes \cite{calderbankkantor}. One important fact about PDSs is that they are reversible, that is, if $d \in D$, then $d^{-1} \in D$. For a subset $B \subset G$, define $B^{(-1)} = \{b^{-1} | b \in B\}$. The symmetric property of PDSs implies $D=D^{(-1)}$. 

We say that a simple, undirected graph $S=(V,E)$ with $v=|V|$ is strongly regular if all vertices have degree $k$ and for any distinct vertices $v_1, v_2 \in V$, $|\{v_3 : (v_1,v_3) \in E \textrm{ and } (v_2,v_3) \in E\}|$ is $\lambda$ if $(v_1,v_2) \in E$ and $\mu$ otherwise. Haemers and Spence (2001) found that there are precisely 167 nonisomorphic SRGs with parameters $(v,k,\lambda,\mu)=(64,18,2,6)$ \cite{HS}. A PDS $D$ can be used to construct an SRG by labeling the vertices of a graph with group elements, and creating an edge between vertices with labels $g,h \in G$ if $gh^{-1} \in D$.

Negative Latin square type (NLST) PDSs are PDSs with the dimensions $(n^2,r(n+1),-n+r^2+3r,r^2+r)$. Note that NLST PDSs with $r=\frac{n-2}{2}$ can be created by taking a reversible Hadamard difference set with the identity and removing the identity element; however, we are not interested in PDSs formed in this way. 

Some constructions of NLST PDSs exist for abelian groups \cite{HIR} \cite{DEN}; for example, Davis and Xiang (2004) found sets of reversible Hadamard difference sets such that the union of translates of these difference sets was a NLST PDS \cite{dx1}, and Polhill, Davis, and Smith (2012) found a NLST PDS in $C_8 \times C_8$  \cite{PDSPDS}. However, very little is known about NLST PDSs in nonabelian groups. Recently, Feng, He, and Chen derived nonabelian PDSs from the Davis-Xiang family \cite{FHC}.

Let $D$ be a NLST PDS. Let $\mathbb{Z}[G]$ be the group ring of $G$ over the integers, that is, the set of all formal sums of elements of $G$. Then in the group ring define $D:=\sum_{d \in D} d$. Notice that we use $D$ to denote both a PDS and the sum of the elements of a PDS in a group ring. This is a standard abuse of notation, and context makes it clear which definition we are referring to. In the group ring, by definition of PDS we have that \begin{equation} 
\label{gringdset}
D^2 = D D^{(-1)} = k 1 + \lambda D  + \mu (G-D-1) \end{equation}

A standard tool for studying PDSs in abelian groups is character theory. A character $\chi$ is a homomorphism that maps an abelian group to the group of complex numbers under multiplication.  To apply characters to elements in the group ring, we extend characters by linearity. For a PDS $D$, we define the character sums of $D$ to be the possible values of $\chi(D)$. The principal character $\chi_0$ maps all elements of a group to 1. Thus $\chi_0(D)=\chi_0(\sum_{i=1}^{k} d_i) = \sum_{i=1}^k \chi_0(d_i) = \sum_{i=1}^k 1 = k$.
 
Now let $\chi$ be a nonprincipal character of a group. One property of nonprincipal characters is that $\chi(G)=0$. Applying $\chi$ to both sides of (\ref{gringdset}) yields $\chi(D)^2 = k \chi(1) + \lambda \chi(D) + \mu \chi(G) - \mu \chi(D) - \mu \chi(1) = (k-\mu) + \chi(D) (\lambda-\mu)$. Using the quadratic formula to solve for $\chi(D)$ yields $\chi(D)= \frac{(\lambda-\mu) \pm \sqrt{(\lambda-\mu)^2 - 4(\mu-k)}}{2}$.  If we restrict our attention to NLST PDSs, then we get that $\chi(D)=r \textrm{ or } r-n$. 

For the rest of the paper we restrict our attention to the following case: let $G$ be a group of order 64 with $k=18$ (so $r=2,n=8$) with an associated normal subgroup $N$ so that the factor group $G/N$ is isomorphic to $C_2^3$. Let $\phi$ be the homomorphism $\phi(g)=Ng$. To apply $\phi$ to the group ring we extend $\phi$ by linearity. Then $\phi(D)^2=\phi(D D^{(-1)}) = k 1_{G/N} + \lambda \phi(D) + \mu (\phi(G)-\phi(D) - 1_{G/N})$. Since $\phi(G)$ is now abelian, we can apply characters. Thus $\chi(\phi(D)^2) = \chi(k 1_{G/N} + \lambda \phi(D) + \mu (\phi(G)-\phi(D) - 1_{G/N}))$, which simplifies to $\chi(\phi(D))^2 =  (-4) \chi(\phi(D)) +  (12) 1_{G/N}$. Like above, solving for $\chi(\phi(D))$ yields $\chi(\phi(D)) = 2 \textrm{ or } -6$. 

Define $C_2^3:=<x^iy^jz^k,x^2=y^2=z^2>$. Let $\chi_i, 1 \le i \le 8$ be the eight possible character mappings on $C_2^3$ (for example, $x\rightarrow 1, y \rightarrow -1, z \rightarrow 1$) and let $f_j$ be an element of $G/N$. Then define an $8 \times 8$ matrix $H:=\{h_{ij} = \chi_i(f_j)\}$ (note that the rows of $H$ are indexed by the characters of $C_2^3$ and the columns are indexed by the elements of $C_2^3$). For example, with the ordering of group elements $C_2^3=\{1,x,y,z,xy,yz,xz,xyz\}$, we have \begin{equation}H= \begin{bmatrix}  
1 & 1 & 1 & 1 & 1 & 1 & 1 & 1  \\
1 & -1 & 1 & 1 & -1 & 1 & -1 & -1 \\
1 & 1 & -1 & 1 & -1 & -1 & 1 & -1 \\
1 & 1 & 1 & -1 & 1 & -1 & -1 & -1 \\
1 & 1 & -1 & -1 & -1 & 1 & -1 & 1 \\
1 & -1 & 1 & -1 & -1 & -1 & 1 & 1 \\
1 & -1 & -1 & 1 & 1 & -1 & -1 & 1 \\
1 & -1 & -1 & -1 & 1 & 1 & 1 & -1 \end{bmatrix} \end{equation}

Let $v$ be a vector of length 8 such that $v_j= |\{d \in D : \phi(d)=f_j\}|$. Then we have that $Hv = [18, \epsilon_2, \epsilon_3, \ldots, \epsilon_8]^T, \epsilon_j \in \{2,-6\}, 2 \le j \le 8$. Taking $H^{-1}$ yields that $v = H^{-1} [18, \epsilon_2, \epsilon_3, \ldots, \epsilon_8]^T, \epsilon_j \in \{2,-6\}, 2 \le j \le 8$.

\end{section}

\begin{section}{Computer Search}

Of the 267 groups of order 64, our algorithm will function for the 212 groups with $C_2^3$ images that are not the elementary abelian group.

To determine if a group $G$ has a $C_2^3$ image, we check if the factor group of any normal subgroup is isomorphic to $C_2^3$, as illustrated by the below algorithm. If we find such a normal subgroup $N$, we store it for later use. 

\begin{algorithmic}

\Procedure{Determine if $G$ has a $C_2^3$ image}{}

\State $G$ = Group \#A group of order 64
\For{normal subgroup $ng$ in $G$}
\If{FactorGroup($ng$)$ \cong C_2^3$}
	\State \#We can run our algorithm on $G$
	\State $N = ng$ \#Store the normal subgroup for later use
	\State break
\EndIf
\EndFor
\EndProcedure
\end{algorithmic}

If $G$ does not have a $C_2^3$ image, we stop the algorithm; otherwise, we continue on. Next, the algorithm calculates all possible values of $v=H^{-1} [18, \epsilon_2, \epsilon_3, \ldots, \epsilon_8]^T, \epsilon_j \in \{2,-6\}, 2 \le j \le 8$, as described by the below procedure. 
\begin{algorithmic}
\Procedure{Find all distributions of PDS elements among the cosets of $G$}{}

\State possiblePDSDistributions = [ ] \#Stores the possible values of $v$

\State Create the character mappings $\chi_i$
\State Let $H=\{h_{ij} = \chi_i(f_j)\}$

\State possibleCharSums = A list of all vectors of length 8 of the form $[18, \epsilon_2, \epsilon_3, \ldots, \epsilon_8]^T, \epsilon_j \in \{2,-6\}, 2 \le j \le 8$
\For{each vector $u$ in possibleCharSums}
	\State Add $H^{-1} u$ to possiblePDSDistributions 
\EndFor
\EndProcedure
\end{algorithmic}

Consider the right cosets of $N$. Note that a value of $v$ corresponds to the distribution of PDS elements among the right cosets of $N$. For example, if $v=[4,2,2,2,2,2,2,2]$, then there are 4 PDS elements in the first right coset of $N$, 2 elements in the second right coset of $N$, so on and so forth. Observe that if $Ng$ is a right coset, and $f \in Ng$, then $f^{-1} \in Ng$ (inverses are always within the same coset). Because inverses are always in the same coset, when we are choosing elements to form a potential PDS we can pick elements out of each coset independently from the other cosets.

We create a database as follows: for each coset, for each number of elements of a hypothetical PDS that belongs in that coset, we compute all of the possible combinations of elements of the coset that satisfy the inverse requirement of PDSs. The creation of this database takes less than a second. We create this database with the below algorithm.

\begin{algorithmic}
\Procedure{Create actualChoicesByCoset array}{}

\State actualChoicesByCoset = []
\For{each coset $C$ of $N$}
\For{$i=0\ldots8$} \#The possible amount of PDS elements that can be contained in $C$
\State actualChoicesByCoset[coset $C$][number of elements $i$] = a list of all of the combinations of size $i$ of PDS elements in coset $C$ that satisfy the inverse requirement
\EndFor
\EndFor
\EndProcedure
\end{algorithmic}

Once we have this database, examining our possiblePDSDistributions, we remove the distributions that are duplicates, contain fractions or negative numbers, or require an order $2$ element in a coset where there is none.

For the search, we choose a vector $v$ from possiblePDSDistributions. We then choose a set from our database's 1st coset section with $v_1$ elements, a set from our database's 2nd coset section with $v_2$ elements, so on and so forth, until we have a set $D$ with 18 elements. As a preliminary check, we check whether $D^2$ contains 18 elements repeated twice and 45 elements repeated six times (and of course 1 element, the identity, repeated 18 times). This check removes most of the PDS candidates. However, this check is not sufficient, as this check does not guarantee that the 18 elements that appear twice in $D^2$ are all elements of $D$. Then, we compute whether (D-2)(D+6)=6G. Any sets that survive this check are indeed PDSs. We considered using a degree of minimal polynomial check, but this check is less time efficient. To speed up these computations we do them in the group ring. We repeat this process for each possible value of $v$. We illustrate this process in the below algorithm. 

\begin{algorithmic}
\Procedure{Search for PDSs in a group $G$}{}
\State pdsList = [] \#Store the PDSs we find
\For{each possible value of $v$}
	\State $temp$ = [ ] \#Store our potential PDS
		\For{each set of elements $a_1$ in the 1st coset of $G$ with size $v_1$}
			\State Add $a_1$ to $temp$
				\For{each set of elements $a_2$ in the 2nd coset with size $v_2$}
					\State Add $a_2$ to $temp$
					\State \vdots
						\For{each set of elements $a_8$ in the 8th coset with size $v_8$}
							\State Add $a_8$ to $temp$
							\If{$temp^2$ contains 45 6s, 18 2s, and 1 18}
								\If{$(temp - 2)(temp+6)=6G$}
									\State Add $temp$ to pdsList \#We have found a PDS!
								\EndIf
							
							\EndIf
							\State Remove $a_8$ from $temp$
						\EndFor

					\State \vdots
					\State Remove $a_2$ from $temp$
				\EndFor
			\State Remove $a_1$ from $temp$

\EndFor
\EndFor
\EndProcedure
\end{algorithmic}

All of the steps except for the last step, ``Search for PDSs in a group $G$", was performed in GAP, a computational group theory programming language \cite{GAP}. GAP was chosen because algebra functions, like groups and images, are built into the language. The last step was performed in Java because Java is significantly faster than GAP, the last step is computationally intensive, and given the precomputation the last step does not require any group theory.

\end{section}

\begin{section}{Results}

\begin{subsection}{Search runtime}

The difficulty with searching for PDSs is the large search space. A true brute force search would require ${64 \choose 18} \approx 3 \times 10^{15}$ checks of possible PDSs. With a PDS check taking roughly $2.2 \times 10^{-5}$ seconds on the author's computer (2019 iMac, 3.6 GHz Intel Core i9, 16 GB RAM) in GAP, searching the entire space of a single group of order 64 would take about 2000 years. We reduce the search space using character theory and inverse checks. 

A typical value of $v$ is $[2,2,2,4,2,0,4,2]$. For all groups except 260, 261, and 266, there are 57 possible values of $v$. Groups 260, 261, and 266, on the other hand, have 92 possible values of $v$ because all of their cosets contain order 2 elements, and so values of $v$ containing only odd numbers are permissible. Thus without inverse checks, it would take $57  {8\choose 2}^5 {8 \choose 4}^2 {8 \choose 0} \approx 5 \times 10^{12}$ checks to fully search a group of order 64. This is notably less than the $3 \times 10^{15}$ that we started with, but would still take around 3 years to complete.

Because PDSs are symmetric, when an element not of order 2 is added to a PDS during a search, then its inverse can also be added ``for free". This significantly reduces the amount of possibilities we need to check. Because of this fact, groups with few order 2 elements take much less time to search than groups with many order 2 elements. The number of order 2 elements serves as an effective proxy when estimating how long a group will take to run. This algorithm would never finish if asked to search the elementary abelian group exhaustively. However, there are already NLST PDSs known in the elementary abelian group \cite{HIR}.

Most groups of order 64 with a $C_2^3$ image have 31 or less order 2 elements. For a rough calculation, we estimate that half of the cosets contain all order 2 elements, and that the other half of the cosets contain no order 2 elements. Thus when elements are chosen from the cosets without order 2 elements, a second element, the inverse, can be added for free. Combined with the character theory search space reduction, it would take about $57 * {8 \choose 2}^3 {8 \choose 4}^1 {8 \choose 0} {4 \choose 1}^2 {4 \choose 2}\approx 1 \times 10^{10}$ PDS checks to fully search a single group, or about 2.5 days.

However, many of these PDS checks are redundant. Consider a set $temp$ with 14 elements (a partially constructed PDS using our algorithm). If $temp^2$ contains any element of the group more than 6 times, we know that $temp$ cannot possibly be a subset of any PDS. Thus we can stop building the PDS early and move on. Stopping the PDS building early saves a lot of time because there might be 16-400 different sets of size 18 that could be built from $temp$. We added checks to determine if the PDS is doomed to fail after adding the elements of the 5th coset, 6th coset, and 7th coset.

Additionally, we implemented the computationally intensive main loop of the program (the actual iterating and building of PDSs) in Java. Combined with the loop optimization, this reduced the time to search all 212 groups exhaustively to about 3 1/2 hours.

\end{subsection}

\begin{subsection}{Search results}

The groups in the GAP library searched for PDSs are as follows: 55-266. All of these groups were searched exhaustively. We found 223,680 PDSs over 49 groups, 1 abelian group and 48 nonabelian groups. Although a NLST PDS was already known in the abelian group 192 ($C_4^2 \times C_2^2$) \cite{dx1}, searching this group provided evidence that our code was working properly. Many PDSs within the same group were equivalent to each other: among the 48 nonabelian groups, there were 176 inequivalent PDSs. Note that PDSs were already known in some nonabelian groups \cite{FHC}. Please find the full PDS database and code at 42ABC/Nonabelian-NLST-PDSs-of-order-64-code-data on Github.

\end{subsection}

\begin{subsection}{Disjoint PDSs}

Using the PDSs generated in 2.1, for each group, we checked if any two NLST PDSs were disjoint, as demonstrated by the below algorithm. NLST PDSs in particular are interesting because in our case, if we have two disjoint NLST PDSs with 18 elements each, the remaining elements in the group are a reversible different set, and so we have a rank 4 Schur ring \cite{BH}. 

\begin{algorithmic}
\Procedure{Find disjoint PDSs}{}
\For{each group $G$ with NLST PDSs}
	\For{each pair of PDSs $D_1,D_2$}	
		\If{$D_1$ and $D_2$ are disjoint}
			\State \#We have found a pair of disjoint PDSs
			\State break
		\EndIf
	\EndFor					
\EndFor
\EndProcedure
\end{algorithmic}

We found that 18 nonabelian and 1 abelian group contain 2 disjoint PDSs. 

\begin{center}
\begin{ssmall}
\begin{tabular}{c|c}
Group & PDSs \\
\hline
60& [ 2, 12, 3, 14, 4, 18, 25, 61, 27, 62, 35, 63, 43, 44, 45, 58, 59, 60 ], [ 5, 6, 7, 20, 30, 52, 15, 56, 39, 40, 46, 48, 49, 50, 34, 54, 23, 64 ]\\
  67& [ 7, 20, 3, 37, 4, 17, 18, 19, 8, 25, 9, 27, 51, 62, 53, 55, 59, 64 ], [ 5, 22, 2, 12, 30, 52, 14, 16, 39, 40, 41, 57, 24, 46, 13, 54, 43, 58  ]\\
  88& [ 2, 32, 3, 14, 4, 19, 25, 47, 27, 49, 13, 34, 53, 54, 55, 63, 45, 64 ], [ 5, 6, 21, 22, 10, 52, 16, 37, 18, 41, 8, 61, 9, 28, 33, 35, 43, 44  ]\\
  90& [ 2, 11, 3, 14, 4, 19, 25, 47, 27, 62, 13, 34, 53, 54, 55, 63, 45, 64 ], [ 6, 20, 10, 52, 15, 36, 38, 56, 17, 39, 40, 41, 24, 48, 49, 50, 43, 60  ]\\
  92& [ 2, 32, 3, 14, 18, 41, 47, 48, 27, 49, 13, 54, 43, 44, 45, 58, 59, 60 ], [ 5, 20, 21, 42, 10, 30, 38, 56, 39, 57, 24, 25, 28, 29, 34, 63, 23, 64 ]\\
  99& [ 2, 11, 3, 14, 4, 19, 26, 61, 27, 62, 55, 63, 23, 43, 44, 58, 59, 60 ], [ 5, 20, 21, 42, 10, 30, 15, 56, 17, 40, 8, 46, 9, 51, 33, 35, 45, 64 ] \\
  123& [ 20, 42, 2, 11, 3, 38, 4, 17, 39, 41, 50, 62, 34, 35, 23, 43, 45, 59 ], [ 5, 21, 10, 30, 14, 56, 18, 19, 40, 57, 8, 26, 46, 61, 29, 51, 33, 63  ]\\
  131& [ 2, 12, 3, 37, 4, 18, 25, 61, 27, 62, 33, 55, 23, 43, 45, 58, 60, 64 ], [ 5, 6, 21, 22, 10, 31, 15, 36, 19, 41, 8, 47, 28, 29, 13, 53, 44, 59 ]\\
  167& [ 6, 22, 2, 52, 4, 17, 19, 40, 8, 46, 27, 51, 13, 33, 34, 53, 59, 60 ], [ 5, 21, 10, 11, 3, 15, 16, 38, 24, 48, 9, 49, 35, 54, 55, 63, 45, 58 ]\\
  174& [ 5, 21, 2, 11, 4, 17, 18, 57, 8, 48, 29, 51, 13, 35, 53, 63, 58, 59 ], [ 20, 42, 10, 30, 3, 14, 16, 37, 19, 39, 40, 41, 24, 61, 27, 49, 44, 45 ]\\
  179& [ 6, 22, 2, 10, 3, 14, 15, 56, 26, 47, 9, 50, 13, 35, 53, 63, 58, 60 ], [ 20, 42, 11, 30, 16, 36, 37, 38, 4, 18, 19, 41, 25, 46, 28, 62, 43, 45 ]\\
  192& [ 5, 20, 3, 16, 4, 17, 19, 57, 8, 48, 9, 28, 50, 62, 34, 55, 58, 59 ], [ 21, 42, 2, 10, 11, 30, 14, 37, 18, 39, 40, 41, 24, 61, 33, 54, 23, 60  ]\\
  193& [ 5, 21, 3, 15, 4, 17, 18, 57, 8, 48, 9, 29, 49, 62, 35, 55, 58, 59 ], [ 20, 42, 2, 10, 12, 31, 14, 36, 19, 39, 40, 41, 24, 61, 33, 53, 23, 60  ]\\
  202& [ 5, 6, 2, 32, 3, 37, 4, 40, 9, 29, 49, 62, 33, 53, 55, 63, 44, 45 ], [ 7, 42, 10, 52, 14, 16, 18, 57, 8, 24, 25, 46, 13, 34, 35, 54, 23, 60 ] \\
  211& [ 5, 42, 3, 15, 4, 19, 8, 26, 27, 49, 13, 53, 54, 55, 44, 59, 60, 64 ], [ 20, 21, 2, 10, 32, 52, 14, 36, 17, 40, 46, 61, 29, 51, 23, 43, 45, 58  ]\\
  226& [ 3, 14, 16, 36, 4, 17, 18, 57, 24, 61, 27, 50, 54, 63, 43, 44, 45, 59 ], [ 5, 21, 2, 10, 30, 32, 37, 56, 8, 25, 9, 29, 49, 62, 34, 35, 23, 60 ]\\
  236& [ 5, 20, 3, 16, 4, 18, 40, 57, 8, 48, 9, 27, 50, 51, 13, 55, 43, 59 ], [ 21, 42, 2, 10, 11, 30, 14, 37, 24, 61, 28, 29, 49, 62, 53, 54, 23, 45  ]\\
  242& [ 14, 37, 4, 17, 18, 40, 8, 24, 25, 61, 9, 51, 34, 53, 54, 55, 45, 60 ], [ 5, 42, 2, 10, 11, 31, 3, 15, 26, 48, 27, 28, 29, 62, 13, 35, 43, 59  ]\\
  261& [ 4, 18, 5, 10, 21, 30, 13, 23, 45, 55, 14, 56, 39, 40, 50, 62, 58, 64 ], [ 2, 11, 3, 8, 26, 38, 9, 19, 28, 41, 20, 42, 24, 47, 17, 57, 53, 54]\\

\end{tabular}
\end{ssmall}
\end{center}

\end{subsection}

\begin{subsection}{Inequivalent PDSs}

For each group, we determined the equivalence classes of the PDSs. For the first PDS $D_1$ in our list, we took all of the automorphisms of that PDS, and found all of the PDSs in our list that were isomorphic to $D_1$. Then, we found the first PDS $D_2$, if any, that was not isomorphic to $D_1$, and repeated this process, as described in the below algorithm. 

\begin{algorithmic}
\Procedure{Find Inequivalent PDSs}{}
\State let $pdsList$ be the list of all the PDSs in this group
\State Give each PDS $D \in pdsList$ a boolean field $D.found$, and set it to false
\For{$D \in pdsList$}
	\If{$D.found$ is false}
		\State Give $D$ a unique ID $i$
		\State Set $D.found$ to true

		\For{each automorphism of $G$}	
			\State Set Aut($D$).$found$ to true
			\State Give Aut($D$) the ID $i$ as well (because $D$ and Aut($D$) are in the same isomorphism class) 
	
		\EndFor
	\EndIf
						
\EndFor
\EndProcedure
\end{algorithmic}

We found that in the abelian group 192, there was only 1 PDS up to equivalence. In the remaining 48 nonabelian groups, there were 176 PDSs up to equivalence. The number of inequivalent PDSs ranged from 1 (group 239) to 13 (group 90), depending on the group.

\end{subsection}

\begin{subsection}{Strongly Regular Graphs}

Using GRAPE \cite{GRAPE}, we determined that our PDSs created eight nonisomorphic SRGs. We compared these isomorphism classes with the SRGs found by Haemers \& Spence (2001) \cite{HS}; all eight of the SRGs we got were in Haemers \& Spence's list. Finding the isomorphism classes and comparing with Haemers \& Spence's list took about 4 minutes. Finding that our PDSs are contained in known (64,18,2,6) SRGs reinforces the validity of our findings and Haemers \& Spence's SRG list. Further, having eight nonisomorphic SRGs suggests that some of these SRGs may be genuinely nonabelian, that is, only producible from a PDS in a nonabelian group. 

\begin{algorithmic}
\Procedure{Find Inequivalent SRGs}{}
\For{each inequivalent PDS}
	\State pds\textunderscore graph = Graph(PDS)
	\For{each SRG in Haemers \& Spence's list}	
		\If{pds\textunderscore graph and SRG are isomorphic}
			\State \#We have found an isomorphic graph
			\State break
		\EndIf
	\EndFor					
\EndFor
\EndProcedure
\end{algorithmic}

\end{subsection}

\begin{subsection}{Hadamard difference set breakdown of PDSs}

Based on some internal structure observed in the PDSs in \cite{dx1}, we searched our database for similar structure that could ultimately be used to generalize these PDSs. Let $R$ be a normal subgroup of $G$ such that $G/R$ is isomorphic to $C_2^2=<z^i w^j | z^2=w^2=1>$. For each inequivalent PDS $D$ from our search, we determined if $D$ could be broken down into three (16,6,2) Hadamard difference sets $A_1, A_2, A_3 \in R$ such that $D=z * A_1 \cup w * A_2 \cup zw * A_3$.

\begin{algorithmic}
\Procedure{Find Difference Set Breakdown}{}
\For{each inequivalent PDS $D$}
	\For{each factor group $G/R$}	
		\If{$G/R$ is isomorphic to $C_2^2=<z^i w^j, z^2=w^2=1>$}
			\State Let $\phi : G \rightarrow G/R$
			\If{$\phi(D)$ contains six copies of $z$, six copies of $w$, and six copies of $z * w$}
				\State \#We separate the elements of $D$ by which coset they belong to
				\State $A_1 = \{d \in D | \phi(d)=z\}, A_2 = \{d \in D | \phi(d)=w\}, A_3 = \{d \in D | \phi(d)=z * w\}$
				\If {the multiset $\{a_1 a_2^{-1} | a_1,a_2 \in A_i\}$ contains every nonidentity element of $R$ twice for $1 \le i \le 3$}
					\State \#We have found a difference set breakdown
					\State break
				\EndIf
			
			\EndIf
		\EndIf
	\EndFor					
\EndFor
\EndProcedure
\end{algorithmic}

We found that, from our list, only PDSs from groups with $C_2^4$ images had a difference set breakdown, and all but 3 of the PDSs from groups with $C_2^4$ images had difference set breakdowns. There was not a clear pattern between SRG isomorphism and whether or not a PDS could be broken down into difference sets. The PDSs that could be broken down into three difference sets spanned several, but not all, SRG isomorphism classes, and some SRG isomorphism classes contained some PDSs that could be broken down and some that could not. 

However, the fact that many of the PDSs can be broken down into 3 difference sets suggests that there may be a construction method for nonabelian NLST PDSs like demonstrated in Davis and Xiang (2004) for some abelian groups \cite{dx1}. This construction method could relate to the work done by Feng, He, and Chen (2020) in nonabelian groups with PDSs with the same SRG as the NLST PDS in $C_4^2 \times C_2^2$ \cite{FHC}.

\end{subsection}

\end{section}

\begin{section}{Conclusions}

In this paper, we found novel NLST PDSs in nonabelian groups of order 64 through a combination of theoretical and computational techniques. 

In future work, the methods recently used by Jedwab and Li to pack PDSs in abelian groups could potentially be adapted and applied to these nonabelian PDSs, aiding in the search for families of NLST PDSs in nonabelian groups \cite{JL}.  

The code used to find the results in this paper could be expanded and improved. A $C_4 \times C_4$ or $C_8 \times C_2$ image version of the code would allow more groups of order 64 not already searched by the $C_2^3$ code to be searched for PDSs. Adjusting the code to search for PDSs with $k=14$ and $k=21$ would be a first step towards creating a database of all PDSs of order 64. 

In experimenting with writing a $C_4 \times C_4$ image code, two interesting objects were found in $G=C_8 \times C_8=<x^iy^j | x^8=y^8=1>$:

\begin{center}
\begin{ssmall}
\begin{tabular}{c|c}
Object & Indices in GAP, Object in group element form \\
\hline
$fakePDS$ &  [ 19, 41, 20, 42, 39, 40, 15, 36, 38, 56, 12, 29, 32, 51, 44, 45, 46, 47]  \\
& $\{x^2 y^4, x^6 y^4, x^4 y^2, x^4 y^6, x^6 y^2, x^2 y^6, x^4 y, x^4 y^3, x^4 y^5, x^4 y^7, x y^4, x^3 y^4, x^5 y^4, x^7 y^4, x^7 y, x^3 y^5, x^5 y^3, x y^7\}$ \\ \hline
$realPDS$ &   [ 12, 15, 17, 19, 20, 23, 24, 29, 32, 36, 38, 41, 42, 51, 56, 57, 60, 61 ]  \\
& $\{x y^4, x^4 y, x^2 y^2, x^2y^4, x^4y^2, x^3 y, x y^3, x^3 y^4, x^5 y^4, x^4 y^3, x^4 y^5, x^6 y^4, x^4 y^6, x^7 y^4, x^4 y^7, x^6 y^6, x^7 y^5, x^5 y^7 \}$ \\
\end{tabular}
\end{ssmall}
\end{center}

$realPDS$ is a NLST PDS with dimensions (64,18,2,6). Note that a NLST PDS was already known in $C_8 \times C_8$ \cite{PDSPDS}. $fakePDS$ is not a PDS, but it somehow has the property that, in the group ring, $realPDS^2=fakePDS^2$. This means that $fakePDS^2$ contains 2 copies of each element in $realPDS$, 18 copies of the identity, and 6 copies of all other elements in $G$. Future work would consist of explaining the existence of $fakePDS$ and finding examples of such objects in other groups.

\end{section}

\begin{section}{Acknowledgements}

The author would like to thank Dr. James Davis for proposing this project, providing advice and suggestions throughout the research, providing the TeX template for the paper, and helping to edit the paper, and Dr. Ken Smith for significant GAP support and theoretical suggestions. The author would also like to thank Drs. Ken Smith and Ryan Kaliszewski for the incidence matrix, symmetric design, and convolution table functions that the algorithms discussed in this paper implement. The author would like to thank Dr. John Polhill for advice on the project, Dr. Prateek Bhakta for recommending loop optimization to improve the runtime of the main PDS search procedure, and Sreya Aluri and Nikita Morozov for suggesting to avoid recursive function calls and to increment an index instead of repeatedly calling a Remove command in the main PDS search procedure. This work was supported by University of Richmond.

\end{section}

\bibliography{NLST_nonabelian_PDS}{}
\bibliographystyle{plain}

\end{document}